\documentclass[12pt]{amsart}
\usepackage{amssymb,amsmath,layout,enumerate}
\newtheorem{theorem}{Theorem}[section]

\newtheorem{prop}[theorem]{Proposition}
\newtheorem{lemma}[theorem]{Lemma}
\newtheorem{cor}[theorem]{Corollary}

\newtheorem{obs}[theorem]{Observation}
\newtheorem{dfn}[theorem]{Definition}

\newtheorem{exm}[theorem]{Example}

\def\proof{\noindent{\bf{Proof.} }}
\def\sqr#1#2{{\vcenter{\hrule height.#2pt
        \hbox{\vrule width.#2pt height#1pt \kern#1pt
                \vrule width.#2pt}
        \hrule height.#2pt}}}

\newcommand{\m}{\mathfrak{m}}
\newcommand{\q}{^{[q]}}
\newcommand{\f}[1]{\ensuremath{\mathfrak{#1}}}
\newcommand{\core}[1]{\ensuremath{{\rm{core}}(#1)}}

\DeclareMathOperator{\depth}{\textrm{depth}}

\begin{document}

\title{A formula for the $*$-core of an ideal}

\thanks{A. N. V. was partly supported by NSA grant H98230-09-1-0057}

\author{Louiza Fouli}
\address{New Mexico State University, Las Cruces, NM, 88003 }
\email{lfouli@math.nmsu.edu}

\author{Janet C. Vassilev}
\address{University of New Mexico, Albuquerque, NM,  87131}
\email{jvassil@math.unm.edu}

\author{Adela N. Vraciu}
\address{University of South Carolina, Columbia, SC, 29208}
\email{vraciu@math.sc.edu}

\subjclass{13A30, 13A35, 13B22}

\keywords{tight closure, reduction, core, *-independent, spread}

\begin{abstract}{Expanding on the work of Fouli and Vassilev \cite{FV}, we
determine a formula for the $*$-$\rm{core}$ of an ideal in two
different settings: (1)  in a Cohen--Macaulay local ring of characteristic $p>0$, perfect residue
field and test ideal of depth at least two, where the ideal has a minimal $*$-reduction that is a 
parameter ideal and (2)  in a normal local domain of characteristic $p>0$, perfect residue
field and $\m$-primary test ideal, where the ideal is a sufficiently high Frobenius power of an ideal. We also exhibit some examples
where our formula fails if our hypotheses are not met.}
\end{abstract}

\maketitle

\section{Introduction}
In a recent paper \cite{FV}, Fouli and Vassilev introduced the
$cl$-$\rm{core}$ of an ideal for a closure operation $cl$.  Like the
core of an ideal originally defined by Rees and Sally \cite{RS}, the
$cl$-$\rm{core}$ of an ideal is the intersection of all the
$cl$-reductions of an ideal, where $J$ is a $cl$-reduction of an
ideal $I$ if $J\subseteq I$ and $J^{cl}=I^{cl}$. Of particular
interest for rings of characteristic $p$, is the $*$-$\rm{core}$,
where $*$ denotes tight closure.

The tight closure of an ideal is contained in the integral closure,
but is generally much closer to the ideal than the integral closure.
Note that a $*$-reduction of an ideal $I$ is also a reduction.
Hence, $\core{I} \subseteq *$-$\core{I}$.  Recall that $\core{I}$
sits deep within $I$; $*$-$\core{I}$ will lie within $I$, but will
be generally much closer to $I$.  In \cite{FV}, Fouli and Vassilev
determined that the $*$-$\rm{core}$ and the core of an ideal $I$
agree if the analytic spread of $I$ is equal to the $*$-spread of $I$, where
the $*$-spread is the minimal number of generators of any minimal
$*$-reduction of $I$.   Fouli and Vassilev also exhibited some
examples where the core and the $*$-core of an ideal are not equal.

Huneke and
Swanson \cite{HS1} determined that in a $2$-dimensional regular local ring the core of an integrally closed ideal $I$  is given by $\core I=I(J:I)$ where $J$ is a minimal
reduction of $I$. This result was later generalized by  Corso, Polini and Ulrich \cite{CPU2} under suitable assumptions on $I$ (not necessarily integrally closed) in the case $(R,\m)$ is a Gorenstein local ring. In particular, if $I$ is an $\m$-primary ideal with reduction number at most $1$ then this formula holds \cite[Theorem~2.6]{CPU2}. We will discuss the additional assumptions on $I$ in Section~\ref{prel}.

Various authors have developed techniques for studying the core of an ideal (see for example \cite{CPU1}, \cite{CPU2}, \cite{HT}, \cite{PU}, \cite{HySm1}, \cite{HySm2}). However, these techniques do not apply when studying the $*$-core and minimal $*$-reductions. For one the analytic spread of an ideal $\ell (I)$ is always bounded above by the dimension of the ring and the $*$-spread can grow arbitrary large.  In order to better understand minimal $*$-reductions and the $*$-core of an ideal we make use of  special tight closure techniques introduced by Vraciu \cite{Vr1} and
further developed in \cite{HV}, \cite{Ep} and \cite{EV}.  To work
with special tight closure, we need to assume that $R$ is an
excellent normal  local ring of characteristic $p>0$ with  perfect residue
field.  We obtain a formula for the $*$-core of an ideal
through special tight closure that is very similar to the formula
shown  by Corso, Polini and Ulrich  in \cite[Theorem~2.6]{CPU2} however with
minimal $*$-reductions instead minimal reductions.  The ideals we
consider in this paper have reduction number equal to one with
respect to their minimal $*$-reductions.  It may be that when
$J\subset I\subset J^*$, the reduction number of $I$ with respect to
$J$ is always one.


The outline of the paper is as follows:  In
Section 2, we briefly introduce the reader to tight closure and conditions associated to the study of the core. In
Section 3 over an excellent normal  local  ring, we provide two
criteria that when  satisfied  we prove that
$*$-$\core I =J(J:I)=I(J:I)$, where $J$ is any minimal $*$-reduction of an ideal $I$ (Theorem~\ref{main}).  We subsequently use these criteria to show that when $(R,\m)$ is a Cohen--Macaulay local ring of characteristic $p>0$
with perfect residue field and test ideal of depth at least two, then $*$-$\core{I}=J(J:I)=I(J:I)$, where $I$ is an ideal that has a minimal $*$-reduction $J$ generated by part of a system of parameters (Theorem~\ref{mainsop}).
If $(R, \m)$ is a normal local ring  of characteristic $p>0$ and test ideal equal to $\m$, we
show that $*$-$\core I=J(J:I)=I(J:I)$ for all ideals $I$ of $R$ (Theorem~\ref{mainm}). Furthermore, when $(R,\m)$ is a normal local ring of characteristic $p>0$, with perfect residue field and $\m$-primary test ideal, then we determine that $*$-$\core{I^{[q]}}=J^{[q]}(I^{[q]}:J^{[q]})$ for $q$ a sufficiently
large power of $p$ and  for every minimal $*$-reduction $J$ of $I$(Theorem~\ref{main*}) .  In Section 4, we exhibit two examples where our formula fails, suggesting that if the criteria of Theorem~\ref{main} are not satisfied then the conclusion need not hold.

\section{ Preliminaries}\label{prel}

Since the $*$-$\rm{core}$ of an ideal $I$ is based on reductions
stemming from the tight closure of $I$, we will review some tight
closure concepts to help clarify these notions for the reader.  For
a more extensive introduction, see \cite{Hu} or \cite{HH1}.

\begin{dfn}{\rm
Let $R$ be a Noetherian ring of characteristic $p>0$. We
denote positive powers of $p$ by $q$ and the set of elements of $R$ which
are not contained in the union of minimal primes by $R^0$. Then
\begin{enumerate}[{\bfseries (a)}]
\item  For any ideal $I \subseteq R$, $I^{[q]}$ is
the ideal generated by the $q$th powers of elements in $I$.

\item We say
an element $x \in R$ is in the {\it tight closure}, $I^{*}$, of $I$ if
there exists a $c \in R^0$, such that $cx^q \in I^{[q]}$ for all large
$q$.

\end{enumerate}}
\end{dfn}

Finding the tight closure of an ideal would be hard without test
elements and test ideals.  A {\it test element} is an element $c \in
R^0$ such that $c I^* \subseteq I$ for all
$I \subseteq R$. Note that $c \in \bigcap\limits_{I \subseteq R}
(I:I^*)$. The ideal $\tau=\bigcap\limits_{I \subseteq R} (I:I^*)$ is
called the {\it test ideal} of $R$.

In a Noetherian local ring of characteristic $p>0$, Vraciu \cite{Vr1}
defined the special tight closure, $I^{*sp}$, to be the elements $x
\in R$ such that $x \in (\frak{m}I^{[q_0]})^*$ for some $q_0=p^{e_0}$.
Huneke and Vraciu show in \cite[Proposition~4.2]{HV} that $I^{*sp}
\cap I=\frak{m}I$ and $I^*=I+I^{*sp}$ if $I$ is generated by
$*$-independent elements. Recall that the elements $f_1, \ldots,
f_n$ are said to be $*$-independent if $f_i \not\in (f_1, \ldots,
\hat{f_i}, \ldots, f_n)^{*}$ for all $i=1, \ldots, n$.

Note that the minimal $*$-reductions of $I$ are generated by
$*$-independent elements. In particular, Epstein shows that if $R$ is an
excellent analytically irreducible local domain then every minimal $*$-reduction
is generated by the same minimal number of generators, namely the
$*$-spread of $I$ (\cite[Theorem~5.1]{Ep}). He also showed in \cite[Lemma~3.4]{Ep}
that $I^{*sp}=J^{*sp}$ for all $*$-reductions of $I$.

We mention here the following result of Aberbach that
is used in many arguments of this paper:

\begin{prop} {\rm (\cite[Proposition~2.4]{Ab})}\label{Ab}
Let $(R, \f{m})$ be an excellent, analytically irreducible local ring of
characteristic $p>0$, let $I$ be an ideal, and let $f\in R$. Assume that
$f\not \in I^*$. Then there exists $q_0=p^{e_0}$ such that for all $q \geq
q_0$ we have $I^{[q]}:f^q \subseteq \f{m}^{[q/q_0]}$.
\end{prop}

\bigskip

We conclude the section by recalling some definitions and conditions associated with the study of the core. We also provide some instances for when these conditions are satisfied.

Let $(R, \m)$ be a Noetherian local ring and $I$ an ideal in $R$. Let $d=\dim R$.

\begin{enumerate}[{\bfseries (a)}]
\item The {\it analytic spread} $\ell=\ell(I)$ is the dimension of the special fiber ring of $I$, $\mathcal{F}(I)= R[It]/\m R[It]=R/\f{m} \oplus I/\f{m}I \oplus I^2/\f{m}I^2 \oplus \ldots$ . It is
well known that $\ell(I)=\mu(J)$ for every minimal reduction $J$ of $I$, where $\mu(J)$ denotes the minimal number of
generators of $J$ \cite{NR}. 

\item $I$ satisfies $G_{s}$ if $\mu(I_p) \leq s-1$ for every prime
$\mathfrak{p} \supseteq I$, with dim$R_{\mathfrak{p}} \leq s-1$. If $I$ is $\m$-primary then $I$ satisfies $G_{d}$ automatically. More generally if $I$ is equimultiple, i.e. $\ell={\rm ht} \; I$, then $I$ satisfies $G_{\ell}$. This condition is rather mild and there are many classes of ideals that would satisfy this condition (see for example \cite{CPU1}, \cite{CPU2}).

\item $K$ is a geometric $s$-residual intersection of an ideal $I$ if there exists an $s$-generated
ideal $\mathfrak{a}$ with $K=\f{a}:I$ and ${\rm ht} \; K \geq s \geq \rm{ht}\; I$ and
${\rm ht}\; I+K \geq s+1$.

Residual intersections are a generalization of linkage theory, when $R$ is a local Gorenstein ring or $I$ is an unmixed ideal. In this case if $s=g$ then $K$ is a $g$-residual intersection of $I$ is equivalent to  $K$ is  linked to $I$ since the ideal $\f{a}$ will be a complete intersection. This topic has been developed in great detail by various authors (see \cite{CEU}, \cite{JU}, \cite{U}).

\item $I$ is  weakly $s$-residually S$_2$ if $R/K$ satisfies Serre's condition S$_2$ for every
geometric $i$-residual intersection with $i \leq s$. Again we remark here that if $I$ is $\m$-primary, or more generally equimultiple then $I$ is weakly $\ell-1$ residually S$_2$.

\end{enumerate}

\section{Criteria and formulas for $*$-core}

The following observation is an immediate consequence of the
definition of the test ideal:

\begin{obs}\label{testideal}
Let $R$ be a ring of characteristic $p>0$, with test ideal $\tau$. Then $\tau I \subseteq *$-$\core{I}$ for every ideal $I\subseteq  R$.
\end{obs}
\proof
Let $J$ be an arbitrary $*$-reduction of $I$. Since $I\subseteq
J^*$, it follows that $\tau \subseteq (J:I)$, and thus $\tau I
\subseteq J$.
\qed

We will use the following characterization of minimal $*$-reductions from \cite{Vr2}:

\begin{theorem}\label{chain}{\rm(\cite[Theorem~1.13]{Vr2})}
Let $(R, \m)$ be an excellent normal local ring of characteristic $p>0$ and perfect residue field.  Let $K\subseteq I$ be
ideals, such that $K$ is generated by $l$ elements, where $l$ is the
$*$-spread of $I$. Assume that $I$ is tightly closed.

The following are equivalent:
\begin{enumerate}[$($a$)$]

\item $K$ is a minimal $*$-reduction of $I$

\item $I^{*sp} \cap K \subseteq \m I$

\item $I=I^{*sp} +K$.
\end{enumerate}

\end{theorem}

Using Theorem~\ref{chain} we can describe explicitly all minimal $*$-reductions of $I$.

\begin{lemma}\label{minred}
Let $(R, \m)$ be an excellent normal local ring of characteristic
$p>0$ and  perfect residue field. Let $I$ be an ideal and let $J$ be a minimal $*$-reduction of $I$. Then we may write  $I=J+(u_1, \ldots, u_s)$, with $u_1, \ldots, u_s \in J^{*sp}$. Suppose that  $f_1, \ldots, f_n$ is a minimal system of generators of $J$ and let $K\subseteq I$. Then $K$ is a minimal $*$-reduction of $I$ if
and only if $K=(f_1+v_1, f_2+v_2, \ldots, f_n+v_n)$, with $v_1,
\ldots, v_n \in (u_1, \ldots, u_s)$.
\end{lemma}
\proof
We check that for every $K=(f_1+v_1, \ldots, f_n+v_n)$ as in the statement,
condition (b) in Theorem~\ref{chain} holds. Indeed, if
$c_1(f_1+v_1)+\ldots + c_n (f_n+v_n)\in I^{*sp}$, then we must have
$c_1f_1 + \ldots + c_n f_n \in I^{*sp}\cap J \subseteq \m J$, where
the last inclusion follows from \cite[ Proposition~1.10]{Vr2}, noting that
$J$ is a minimal $*$-reduction of $I$.

Conversely, let $K\subseteq I$ be a minimal $*$-reduction.
Since $J$ is a minimal $*$-reduction of $I$, we know that the $*$-spread
of $I$ is $n$, the minimal number of generators of $J$. Thus, $K$ must
also be minimally generated by $n$ elements, and we can write
$K=(g_1+v_1, \ldots, g_n + v_n)$, with $g_i \in (f_1, \ldots, f_n)$, and
$v_i\in (u_1, \ldots, u_s)$ for all $i=1, \ldots, n$.

We claim that we must have $(g_1, \ldots, g_n)=(f_1, \ldots, f_n)$.
Assuming this claim, we can write $\underline{f}=A\underline{g}$,
with $A$ an invertible $n \times n$ matrix with entries in $R$,
where $\underline{f}$ denotes the vector $(f_1, \ldots, f_n)^{\rm T}$. The
entries of $A(g_1+v_1, \ldots, g_n+v_n)^{\rm T}$ are then generators for
$K$, and they can be written in the form $f_1+v_1', \ldots, f_n
+v_n'$, where $\underline{v'}=A\underline v$.

Assuming by contradiction that the claim is not true, we can choose
a minimal system of generators $f_1', \ldots, f_n'$ for $J$ such
that $(g_1, \ldots, g_n)\subseteq \m J +(f_2', \ldots, f_n')$. In
order to achieve this, one can take $f_1', \ldots, f_n'$ to be the
pull-backs to $R$ of a vector space basis $\overline{f_1'}, \ldots,
\overline{f_n'}$ of $J/\m J$ obtained by extending a basis of the
subspace spanned by $\overline{g_1}, \ldots, \overline{g_n}$.

Since $K$ is a minimal $*$-reduction of $I$, we must have $f_1'\in
K^*$, and therefore $f_1'\in (\m f_1', f_2', \ldots, f_n', u_1,
\ldots, u_s)$. Using Proposition~\ref{Ab}, this implies that
$f_1'\in (f_2', \ldots, f_n', u_1, \ldots, u_s)^*$.  Since $u_1,
\ldots, u_s \in J^{*sp}$, we now see that in fact we must have
$f_1'\in (f_2',\ldots, f_n')^*$. Indeed, there exists $c \in R^{0}$ such that $$cf_1'^q=a_2f_2'^q+\ldots +a_nf_n'^q+b_1u_1^q+\ldots +
b_su_s^q ,$$ for some $a_i, b_j \in R$ with $2 \leq i \leq n$ and $1 \leq j \leq s$. Since $u_j \in J^{*sp}$ for all $j$ there exists $q_0$ such that we can rewrite to above equation to obtain:
$$c^2f_1'^q=ca_2f_2'^q+\ldots +ca_nf_n'^q+cb_1u_1^q+\ldots +
cb_su_s^q =a_1'f_1'^q+ (ca_2+a_2')f_2'^q+\ldots + (ca_n+a_n')f_n'^q,$$
with $a_1', \ldots, a_n'\in \m^{q/q_0}$. The conclusion follows by applying
Proposition~\ref{Ab}. This contradicts the fact that $f_1', \ldots,
f_n'$ are $*$-independent. \qed

\begin{theorem}\label{main}  Let $(R, \m )$ be a normal local ring of
characteristic $p>0$ and perfect residue field. Let $I$ be an ideal and $J$ be a minimal $*$-reduction of $I$. 
Then we may write  $I=J+(u_1, \ldots, u_s)$, with $u_1, \ldots, u_s \in J^{*sp}$. Suppose that  $f_1, \ldots, f_n$ is a minimal system of generators of $J$.

\begin{enumerate}[$($a$)$]
\item If $(f_1, \ldots, \hat{f}_{i}, \ldots, f_n):f_i
\subseteq (J:I)$, and $u_j (J:u_j)\subseteq J(J:I)$ for all $j=1,
\ldots ,s$, then $*$-$\core{I}\subseteq J(J:I)$.

\item If $u_j(J:I)\subseteq \m J(J:I)$ for all $j=1, \ldots, s$, then
$*$-$\core{I}\supseteq J(J:I)$.

\item If the assumptions of both part~$($a$)$ and part~$($b$)$ are satisfied  then
$$*\text{-}\core{I}=J(J:I)=I(J:I).$$

\end{enumerate}

\end{theorem}

\proof (a) Let $A_1f_1+\ldots + A_n f_n \in *$-$\core{I}$. Fix $i \in \{1,
\ldots, n\}$ and $j \in \{1, \ldots, s\}$. We wish to show that $A_iu_j
\in J$. Since $j$ is arbitrary, it will follow that $A_i \in (J:I)$.

Let $J'=(f_1, \ldots, f_{i-1}, f_i+u_j, f_{i+1}, \ldots, f_n)$ be
another minimal $*$-reduction of $I$. We must have $A_1f_1+\ldots
+A_nf_n = B_1f_1 + \ldots +B_{i-1}f_{i-1} + B_i (f_i + u_j) + \ldots
+ B_n f_n$ for some $B_1, \ldots, B_n \in R$. Note we must have $B_i
u_j \in u_j(J:u_j)\subseteq J(J:I)$ by assumption, and thus we can
write $B_iu_j=C_1f_1 + \ldots + C_nf_n$ with $C_1, \ldots, C_n \in
(J:I)$. Therefore  $A_i -B_i
-C_i \in (f_1, \ldots , \hat{f}_{i}, \ldots, f_n):f_i \subseteq
(J:I)$ and since $C_i \in (J:I)$, then $A_i-B_i \in (J:I)$. Multiplying by $u_j$ we obtain $A_iu_j-B_iu_j \in I(J:I)$ and thus $A_iu_j \in I(J:I) \subseteq J$.

(b) Let  $A_1f_1+ \ldots + A_n f_n \in *$-$\core{I}$, with  $A_1, \ldots, A_n \in (J:I)$, and let $J'$ be another minimal $*$-reduction of $I$. Then by Lemma~\ref{minred} we can assume that
$J'=(f_1+v_1, \ldots, f_n+v_n)$, where $v_1, \ldots, v_n \in (u_1,
\ldots, u_s)$. We wish to find $B_1, \ldots, B_n \in (J:I)$ such
that
\begin{equation}\label{eq}
A_1f_1+ \ldots + A_n f_n=B_1(f_1+v_1) + \ldots + B_n (f_n+v_n).\end{equation}

Write $(J:I)=(g_1, \ldots, g_l)$, and $A_i=\Sigma_{j=1}^l\alpha
_{ij}g_j$, $B_i=\Sigma_{j=1}^l y_{ij}g_j$, where $y_{ij}$ are
unknowns. By assumption, we can write $v_ig_k=\Sigma_{\nu=1}^n
\Sigma_{j=1}^l \gamma_{ikj\nu} g_jf_{\nu}$, with $\gamma_{ikj\nu}\in
\m$ for all $i=1, \ldots, n$ and all $k=1, \ldots, l$.

A sufficient condition for equation~(\ref{eq}) to hold is that the coefficients of $f_{\nu}$ on each side
of the equation are the same for all $\nu=1, \ldots, n$, which leads to the system of linear equations
\begin{equation}\label{eq2}
\Sigma_{j=1}^l y_{\nu j} g_j+ \Sigma_{j=1}^l\Sigma_{i=1}^n\Sigma_{k=1}^l y_{ik} \gamma_{ikj\nu} g_j = \Sigma_{j=1}^l \alpha _{\nu j} g_j\ \ \ \ \mathrm{for \ all} \ \nu =1, \ldots, n.
\end{equation}
A further sufficient condition is obtained by identifying the coefficients of each $g_j$
on the two sides of the equations~(\ref{eq2}). The following system of linear equations is thus obtained:
\begin{equation}\label{eq3}
(\nu, j): \ \  \ \ y_{\nu j} + \Sigma_{i=1}^n\Sigma_{k=1}^l y_{ik}\gamma_{ikj\nu}=\alpha _{\nu j} \ \ \ \mathrm{for \ all} \ \nu =1, \ldots, n, \  j=1, \ldots, l.
\end{equation}
This is a system of $nl$ equations with $nl$ unknowns. Note that each
$y_{\nu j}$ appears with a unit coefficient in equation $(\nu, j)$ and with coefficient in $\m$ in all the other equations, and
therefore the matrix of coefficients can be row-reduced to the identity
matrix, showing that there exists a unique solution for the system. \qed

\begin{cor}
With the same set up as in Theorem~\ref{main} if the assumptions of both part~$($a$)$ and part~$($b$)$ are satisfied then  $*$-$\core{I}$ is obtained as a finite
intersection of minimal $*$-reductions, namely
\begin{center}
$*$-$\core{I}=J \cap\left( {\bigcap_{i, j} J_{i, j}}\right),$
\end{center}
where $J=(f_1, \ldots, f_n)$ and $J_{i, j}:=(f_1, ..., f_{i-1},f_i
+u_j, f_{i+1}, ...f_n)$, for all $i\in \{1, \ldots, s\}$ and $j\in \{1, \ldots, s\}$.
\end{cor}
\proof
The proof of part $($a$)$ of Theorem~\ref{main} shows that $J \cap
(\bigcap_{i, j} J_{i, j})\subseteq J(J:I)$, while part $($b$)$ of Theorem~\ref{main} shows that $J(J:I)\subseteq *$-$\core{I}$. The
conclusion follows, since $*$-$\core{I}\subseteq J\cap (\bigcap_{i, j}
J_{i, j})$ by the definition of the $*$-$\rm{core}$.
\qed

\begin{theorem} \label{mainsop}
Let $(R,\m)$ be Cohen Macaulay normal local ring of characteristic $p>0$ and perfect residue field. Let $x_1, \ldots, x_d$ be part of a
system of parameters, with $x_1, x_2\in \tau$, where $\tau$ is the
test ideal of $R$. Let $I$ be an $R$-ideal and suppose $J=(x_1^t, x_2^t, x_3\ldots,
x_d)\subseteq I \subseteq J^*$, where $t \geq 2$. Then
$*$-$\core{I}=J(J:I)=I(J:I)$.
\end{theorem}
\proof First notice that the normality assumption implies that
$\depth{\tau} \geq 2$, and thus we can choose $x_1, x_2 \in \tau$.
This is because the defining ideal of the non-regular locus has
height at least two, and \cite[Theorem~6.2]{HH2} shows that we can
pick a regular sequence consisting of two elements $c, d$ in this
ideal, and replace them by sufficiently large powers in order to
obtain a regular sequence of length two consisting of test elements.

Let $I=J+(u_1, \ldots, u_s)$, where $u_i \in I^{*sp}$ for all $1\leq i \leq s$. We need to  see that all the
hypotheses listed in part~$($a$)$ and part~$($b$)$ of Theorem~\ref{main} hold.

To confirm that the hypotheses of part~$($a$)$ hold notice that

 $$(x_1^t, x_2^t, x_3, \ldots, \widehat{x_{i}}, \ldots, x_d):x_i = (x_1^t, x_2^t, x_3, \ldots,
\widehat{x_{i}}, \ldots, x_d)\subseteq J \subseteq (J:I).$$

For the second hypothesis in part~$($a$)$, note that $$J^*\subseteq J:\tau
\subseteq J:(x_1, x_2)=(x_1^t, x_2^t, x_3,  \ldots, x_d, (x_1x_2)^{t-1}).$$ Thus,
$J^{*sp}\subseteq \m J +(x_1x_2)^{t-1}$, and we may assume without loss of
generality that we can write $u_j=(x_1x_2)^{t-1}u_j'$. We have
$$u_j(J:u_j)=(x_1x_2)^{t-1}u_j'( (x_1, x_2, x_3, \ldots, x_d):u'_j)\subseteq
(x_1, x_2)J \subseteq (J:I)J.
$$
 The last inclusion follows from the fact that $(x_1, x_2)\subseteq \tau \subseteq (J:I)$.

Finally we establish that the hypothesis in part~$($b$)$ holds. As above, we can write $u_j=(x_1x_2)^{t-1} u_j'$, and thus $(J:I)\subseteq (x_1, x_2, x_3, \ldots, x_d):(u_1', \ldots,
u_s')$. Therefore
$$u_j(J:I)\subseteq (x_1x_2)^{t-1}(x_1, x_2, x_3, \ldots, x_d)\subseteq (x_1^{t-1}, x_2^{t-1})J \subseteq \m J (J:I),
$$ where the last inclusion follows from the fact that $(x_1^{t-1}, x_2^{t-1})\subseteq \m (x_1, x_2) \subseteq \m (J:I)$
when we take $t \ge 2$.
 \qed

The following observation can be obtained as a corollary of Theorem~\ref{mainsop} and \cite[Theorem~4.8]{FV},
using \cite[Theorem~2.6 ]{CPU2}. Nonetheless, we give an elementary direct proof.

\begin{obs}\label{red num}
Let $(R, \m)$ be a Cohen-Macaulay normal local ring of characteristic $p>0$,
perfect residue field and  test ideal $\tau$. Let $x_1,\ldots, x_d$
be  part of a system of parameters with $x_1, x_2 \in \tau$. Let $I$
be an $R$-ideal and suppose $J=(x_1^t, x_2^t, x_3,  \ldots, x_d)\subseteq I
\subseteq J^*$, where $t\ge 2$. Then $r_J(I)=1$.
\end{obs}
\proof
As seen in the proof of Theorem~\ref{mainsop}, we can write $I=J+(x_1x_2)^{t-1}(u_1', \ldots, u_s')$.
Thus, it suffices to see that
$(x_1x_2)^{t-1}u_i' \in J^2$ for all $i, j \in \{1, \ldots, s\}$.
This is immediate, since $(x_1x_2)^{2t-2}\in (x_1^tx_2^t) \subseteq J^2$.
\qed

\begin{cor}
Let $(R,\m)$ be a Cohen-Macaulay normal local ring of characteristic $p>0$, perfect
residue field and test ideal $\tau$ . Suppose $x_1,\ldots,
x_d$ form a system of parameters for $R$ and that $x_1, x_2 \in \tau$.  Let $I$ be an $R$-ideal and suppose
$(x_1^t,\ldots,x_d^t) \subseteq I \subseteq (x_1^t, x_2^t, x_3, \ldots, x_d)^*$, where $t \geq 2$. Then
$\text{\rm gr}_{I}R$ is Cohen Macaulay.
\end{cor}

\proof
Notice that with the above conditions $I$ is an $\m$-primary ideal and satisfies $G_{\ell}$, where
$\ell$ is the analytic spread. Also the condition $\depth R/I^{j} \geq d-g-j+1$ for all
$1 \leq j \leq \ell-g +1$, where $g= {\rm ht} (I)$ is immediately satisfied and by Observation~\ref{red num} $r(I)=1$. Therefore, by \cite[Theorem~3.1]{JU} we have ${\rm gr}_{I}(R)$ is Cohen-Macaulay. \qed

\begin{theorem}\label{mainm}
Let $(R,\m )$ be a normal local ring of characteristic $p>0$ and perfect
residue field.  If the test ideal is $\m$, then
$*$-$\core{I}=J(J:I)=I(J:I)$ for any minimal $*$-reduction $J$ of $I$.
\end{theorem}

\proof Since $\tau \subseteq (J:I)$, we have $(J:I)=\m$ or $(J:I)=R$. Note
that $(J:I)=R$ is equivalent to $I=J$, which means that $I$ is the
only minimal $*$-reduction of $I$, and thus the conclusion holds in this case.

Hence, we may assume that $(J:I)=\m$. We have $(J:I)J\subseteq
*$-$\core{I}$ by Observation~\ref{testideal}. In order to
verify the other inclusion, we will check that the assumptions in
part (a) of Theorem~\ref{main} hold. The fact that $(f_1, \ldots,
\hat{f_i}, \ldots, f_n):f_i \subseteq \m$ is immediate, since $f_1,
\ldots, f_n$ is a minimal system of generators of $J$. It remains to
check that $u(J:I)\subseteq \m J$, for $u \in J^{*sp}$. 

Since $u \in J^{*sp} \cap I$ then $u(J:I) \subseteq J^{*sp}$. On the other hand, $u (J:I) \subseteq J$ and thus $u (J:I) \subseteq J^{*sp} \cap J=\m J$, where the last equality holds by \cite[Proposition~4.2]{Vr1}.

\qed

\begin{lemma}\label{special}
Let $(R, \m)$ be an excellent, analytically irreducible local ring of characteristic $p>0$ and $\m$-primary test ideal $\tau$.
Let $J=(f_1, \ldots, f_n)$ and  consider $I=(f_1, \ldots , f_n, u_1, \ldots, u_s)$,
where $u_1, \ldots, u_s \in J^{*sp}$. Then there exists $q_0$ sufficiently
large such that the following hold:
\begin{enumerate}[{\mdseries $($a$)$}]
\item $(f_1^{q}, \ldots, \hat{f_i^q}, \ldots,
f_n^{q}):f_i^{q} \subseteq \tau;$ for all $i=1, \ldots, n$, and
\item $u_k^{q}(J^{[q]}:u_k^{q})\subseteq \m \tau J^{[q]}$ for all $k=1, \ldots, s$
\end{enumerate}
\end{lemma}
\proof Part $($a$)$ follows from Proposition~\ref{Ab}.

It remains to  prove part~$($b$)$. Fix $j$ sufficiently large so that $(\tau ^j :
\tau) \subseteq \m \tau $ (this is possible by the Artin-Rees Lemma,
since $(\tau ^j:\tau) \subseteq \tau ^j:d$ for a fixed non-zerodivisor
$d\in \tau$, and $Fd \in \tau ^j $ implies that $ Fd \in
d\tau^{j-c}$, i.e. $F\in \tau ^{j-h}$, where $h$ is a constant),
and fix $q$ sufficiently large so that $(f_1^{q}, \ldots,
\hat{f_i^q}, \ldots, f_n^{q}):f_i^{q}\subseteq \tau^j$
for all $i=1, \ldots, n$.

Since $u_k\in J^{*sp}$, we can choose $q \gg 0$ so that $u_k^{q}\in
(\tau^j J^{[q]})^*\subseteq (\tau^j J^{[q]}):\tau $. Let $g \in
(J^{[q]}:u_k^{q})$, and write $gu_k^{q}=c_1f_1^{q}+\ldots + c_n
f_n^{q}$. Multiplying by $t \in \tau $ yields $tg
u_k^{q}=(tc_1)f_1^{q} + \ldots + (tc_n)f_n^{q}$. But on the other
hand we have $tu_k^{q}= d_1f_1^{q}+\ldots + d_n f_n^{q}$, with $d_1,
\ldots, d_n \in \tau ^j$. It follows that $tc_i - gd_i\in (f_1^{q},
\ldots, \hat{f_i^q}, \ldots, f_n^{q}):f_i^{q}\subseteq
\tau ^j$, and therefore $c_i\in (\tau ^j : \tau) \subseteq \m \tau $.
\qed

\begin{theorem} \label{main*}
Let $(R, \m)$ be a  normal local ring of characteristic $p>0$, perfect residue field and $\m$-primary test ideal
$\tau $. Let $I$ be an $R$-ideal and let $J \subseteq I$ be a minimal $*$-reduction of
$I$. Then $*$-$\core{I\q}=J\q (J\q :I\q)$ for all $q=p^e \gg 0$.
\end{theorem}
\proof Let $I=J+(u_1, \ldots, u_s)$.  We need to prove $*$-$\core{I\q} =
J\q(J\q:I\q)$ by showing the hypotheses of Theorem~\ref{main} hold when
$I, J$ are replaced by $I\q, J\q$ respectively.

First note that by Lemma~\ref{special}~part~$($a$)$
$(f_1^q,\ldots,\hat{f_i^q},\ldots,f_n^q):f_i^q \subseteq \tau \subseteq
(J:I)$
 showing the first hypothesis of part (a) in Theorem~\ref{main} is satisfied. Using Lemma~\ref{special}~part~$($b$)$ we have that $u_j^{q}(J^{[q]}:u_j^{q})\subseteq \m \tau
J^{[q]} \subseteq J\q(J\q:I\q)$ for all $j=1, \ldots,
 s$ which shows that the hypotheses of Theorem~\ref{main}~part$($a$)$ hold.

In order to verify that the hypothesis of Theorem~\ref{main}~part~$($b$)$ holds we note that $\tau \subseteq (J\q : I\q)$. Thus the result follows by Lemma~\ref{special}~part~$($b$)$.

\qed

\section{Examples}

The first example was an inspiration for Theorem~\ref{mainsop} although
the ring is not normal.

\begin{exm}
{\rm  Let $R=k[[x,y,z]]/(xyz)$, where $k$ is a field of characteristic $p>0$ and $I=(x,yz)$.  By \cite[Theorem~3.7]{Va},
$\tau=(xy,xz,yz)$ and $I=(x+ayz)^*$ for any $0 \neq a \in k$. Note, that
$((x+ayz):I)=(x^2,y^2z^2)= \cap_{a \in k} (x+ayz)=*$-$\core{I}$.}
\end{exm}

If the parameter ideal is not sufficiently embedded for Theorem~\ref{mainsop} to apply, the following example shows that $*$-$\core
I \neq J(J:I)$.

\begin{exm} {\rm Let $R=k[x, y, z]/(x^{5}+y^{5}+z^{5})$, where $k$ is a perfect field of characteristic $p>0$ such that $p \neq 5$. Let $J=(y, z)$ and set $u=x^2$ and $I=J+(u)$.
Note that $u \in J^*$ as can be seen from the fact that $\tau=\m^3$ and
$u\in (J:\m^3)$.  The assumptions in both part~$($a$)$ and part~$($b$)$ of
Theorem~\ref{main} fail, and $J(J:I) \nsubseteq (y+x^2, z)$ which is
a minimal $*$-reduction of $I$.
 }\end{exm}

We are going to use the following result of Brenner to give another
example where the equality $*$-$\core{I}=J(J:I)$ does not always obtain
when the assumptions of Theorem~\ref{main} do not hold.

\begin{theorem}\label{degree}{\rm (\cite[Example~5.2]{Br})}
If $K$ is an algebraically closed field, and
$$R=\frac{K[x, y, z]}{(x^k+y^k+z^k)},$$
where $k$ is not divisible by the characteristic of $K$, and
$f=x^{d_1}, g=y^{d_2}, h=z^{d_3}$ are such that $d_i\le k$ for $i\in
\{1, 2, 3\}$, and $d_1+d_2+d_3\le 2k$, then $R_{\ge k}\subseteq (f,
g, h)^*$.
\end{theorem}

Note that if $R$ is a graded ring, $J$ a homogeneous ideal, and
$u\in J^*$ has degree larger than the degrees of the generators of
$J$, then $u\in J^{*sp}$.

\begin{exm} {\rm Let $R=k[x, y, z]/(x^{10}+y^{10}+z^{10})$, where $k$ is an algebraically closed field of characteristic $p>0$ such that $p \neq 2, 5$. Let $J=(x^5, y^7, z^8)$ and $u=xy^3z^6$.  By Theorem~\ref{degree},
$u$ is in the tight closure of $J$.  Define $I=J+(u)$, and note that
$J$ is a minimal $*$-reduction of $I$. The assumption that
$u(J:I)\subseteq J(J:I)$ fails in Theorem~\ref{main} because
$(J:I)=(x^4, y^4, z^2)$, and $z^2u\notin J(J:I)$. We have that
$J(J:I) \not\subseteq*$-$\core{I}$. More specifically,
$x^5z^2 \in J(J:I)$, but it is not in $J^{\prime}=(x^5+xu, y^7+yu,
z^8+zu)$, which is a minimal $*$-reduction of $I$.
 }\end{exm}


\begin{thebibliography}{MMMM}

\bibitem[Ab]{Ab}{I. Aberbach, {\em Extensions of weakly and strongly F-rational rings by flat maps},
J. Algebra, {\bf  241} (2001), 799--807.}

\bibitem[Br]{Br} Brenner, H., {\em Computing the tight closure in dimension two}, Math. Comp., {\bf 74}  (2005),  no. 251, 1495--1518

\bibitem[CEU]{CEU}{ Chardin,~M., Eisenbud,~D. and Ulrich,~B.,  {\em Hilbert functions, residual intersections, and residually ${\rm S}_2$ ideals}, Compositio Math., {\bf 125}  (2001),  no. 2, 193--219.}


\bibitem[CPU1]{CPU1}
Corso,~A., Polini,~C., Ulrich,~B., {\em The structure of the core of
ideals},  Math. Ann.  {\bf 321}  (2001),  no. 1, 89--105.

\bibitem[CPU2]{CPU2}
Corso,~A., Polini,~C., Ulrich,~B., {\em  Core and residual intersections
of ideals},  Trans. Amer. Math. Soc.,  {\bf 354} (2002), no. 7, 2579--2594.


\bibitem[Ep]{Ep} Epstein,~N., {\em A tight closure analogue of
analytic spread}, Math. Proc. Camb. Phil. Soc., {\bf 139}, (2005),
371--383.

\bibitem[EV]{EV} Epstein,~N., Vraciu,~A., {\em A Length
Characterization of $*$-spread}, Osaka J. Math., {\bf 45}  (2008), no. 2,
445--456.


\bibitem[FV]{FV}
Fouli,~L., Vassilev,~J., {\em The $cl$-core of an ideal}, preprint, arXiv: 0810.3033 [math.AC].

\bibitem[HH1]{HH1}
Hochster,~M., Huneke,~C., {\em Tight closure, invariant theory, and the
Briançon-Skoda theorem}, J. Amer. Math. Soc., {\bf 3} (1990), no. 1,
31--116.

\bibitem[HH2]{HH2} Hochster,~M., Huneke,~C., {\em F-regularity, test elements, and smooth base change}, Trans. Amer. Math. Soc., {\bf 346} (1994), 1--62.

\bibitem[Hu]{Hu}
Huneke,~C., {\em Tight Closure and Its Applications}, CBMS Lect. Notes
Math., {\bf 88}, American Math. Soc., Providence 1996.

\bibitem[HS1]{HS1}
Huneke,~C., Swanson,~I., {\em Cores of ideals in 2-dimensional
regular local rings}, Michigan Math. J., {\bf 42} (1995), 193--208.

\bibitem[HT]{HT}
Huneke,~C. and Trung,~N., {\em On the core of ideals}, Compos. Math., {\bf 141}  (2005), 1--18.


\bibitem[HV]{HV}
Huneke,~C., Vraciu,~A., {\em Special tight closure}, Nagoya Math. J., {\bf
170} (2003), 175--183.

\bibitem[HySm1]{HySm1}{ Hyry,~E. and  Smith,~K {\em On a non-vanishing
conjecture of Kawamata and the core of an ideal}, Amer. J. Math.,
{\bf 125} (2003), 1349--1410.}

\bibitem[HySm2]{HySm2}{Hyry,~E. and  Smith,~K, {\em Core versus graded core, and
global sections of line bundles}, Trans. Amer. Math. Soc. ,{\bf 356} (2003),
3143--3166.}


\bibitem[JU]{JU} Johnson, M., Ulrich, B., {\em Artin-Nagata properties and Cohen-Macaulay associated graded rings},  Compositio Math.,  {\bf 103 } (1996),  no. 1, 7--29.


\bibitem[NR]{NR}{Northcott,~D.~G. and Rees~ D., {\em Reductions of ideals in
local rings}, Proc. Camb. Phil. Soc. {\bf 50} (1954), 145--158.}


\bibitem[PU]{PU}
Polini,~C., Ulrich,~B., {\em A formula for the core of an ideal}, Math.
Ann., {\bf 331}  (2005),  no. 3, 487--503.

\bibitem[RS]{RS}
Rees,~D., Sally,~J., {\em General elements and joint reductions}, Michigan
Math. J.,  {\bf 35} (1988), no. 2, 241--254.


\bibitem[U]{U}{Ulrich,~B., {\em Artin--Nagata properties and reductions
of ideals}, Contemp. Math., {\bf 159} (1994), 373--400.}


\bibitem[Va]{Va}
Vassilev,~J., {\em Test Ideals in quotients of $F$-finite regular local
rings}, Trans. A.M.S., { \bf 350} (1998), 4041--4051.


\bibitem[Vr1]{Vr1}
Vraciu,~A., {\em $*$-independence and special tight closure} , J. Algebra,
{\bf 249}  (2002),  no. 2, 544--565.

\bibitem[Vr2]{Vr2}
Vraciu,~A., {\em Chains and families of tightly closed ideals}, Bull.
London Math. Soc., {\bf 38}  (2006),  no. 2, 201--208.

\end{thebibliography}
\end{document}